\documentclass[11pt]{article}
\usepackage{amsfonts}
\newcommand{\halmos}{\rule{1ex}{1.4ex}}
\newcommand{\proofbox}{\hspace*{\fill}\mbox{$\halmos$}}

\newenvironment{proof}{\noindent {\bf Proof}.}
{\proofbox\par\smallskip\par}

\newcommand{\p}{{\mathbb{P}}}

\newcommand{\np}{{\cal G}_{n,p}}

\newcommand{\Nat}{\mathbb{N}}
\newcommand{\ex}{{\mathbb{E}}}

\newcommand{\eps}{\varepsilon}
\newcommand{\la}{\lambda}

\newcommand{\zi}{{\mathbb{Z}}}

\newcommand{\D}{{\bf D}}

\newcommand{\dd}{{\bf d}}

\newtheorem{theorem}{Theorem}[section]
\newtheorem{lemma}[theorem]{Lemma}

\title{Percolation on sparse random graphs with given degree sequence}

\author{N. Fountoulakis\thanks {The author is supported by the EPSRC, grant no.~EP/D50564X/1} \\
School of Mathematics \\ University of Birmingham \\
Edgbaston, B15 2TT \\ United Kingdom \\ \texttt{nikolaos@maths.bham.ac.uk}}

\begin{document}
\maketitle
\begin{abstract}
We study the two most common types of percolation process on a
sparse random graph with a given degree sequence. Namely, we examine
first a bond percolation process where the edges of the graph are
retained with probability $p$ and afterwards we focus on site
percolation where the vertices are retained with probability $p$. We
establish critical values for $p$ above which a giant component
emerges in both cases. Moreover, we show that in fact these
coincide. As a special case, our results apply to power law random 
graphs. We obtain rigorous proofs for formulas derived by several 
physicists for such graphs. 
\end{abstract}

\section{Introduction}
Traditionally percolation theory has been the study of the properties of
a random subgraph of an infinite graph, that is obtained
by deleting each edge of the graph with probability $1-p$ for some
$p\in (0,1)$ independently of every other edge. The question
that has been mainly investigated is whether the subgraph that is spanned by
these edges has an infinite component or not. The classical type of graphs that
was studied in percolation theory is the lattice
$\mathbb{Z}^d$ in various dimensions $d\geq 2$ (see~\cite{Grim}). Until
now various other types of lattices have been studied. In each of the above
cases the main problem is the calculation of a critical $p_c$ so that
if $p < p_c$ then the random subgraph obtained as above has no
infinite components, whereas if $p> p_c$ there is an infinite
component with probability 1.

In the present work, we study percolation on finite graphs whose
number of vertices is large. This problem is old,
in the sense that for example a $\np$ random graph is a
random subgraph of the complete graph on $n$ vertices, where each
edge appears with probability $p$ independently of every other
edge. In this context, a question about the appearance of an infinite
component is senseless. A somehow analogous question is whether there
exists a component of the random subgraph containing a certain
proportion of the vertices or as we customarily say a {\em giant
component}. More specifically, if the original graph has $n$ vertices
the question now is whether there exists an $\varepsilon >0$ for which
there is a component of the random subgraph that has at least
$\varepsilon n$ vertices with probability $1-o(1)$ (as $n\rightarrow
\infty$). Hence, we also ask (quite informally)
for the existence of a critical $p_c$ for which whenever $p< (1-\delta)p_c$
then for every $\varepsilon >0$ there is no component having at
least $\varepsilon n$ vertices with probability $1-o(1)$ and whenever
$p > p_c(1+\delta)$ then there is a component with at least $\varepsilon n$
vertices for some $\varepsilon >0$ with probability $1-o(1)$.
A classical example of this is the $\np$ random subgraph of $K_n$, 
the complete graph on $n$ vertices, where
as it was proved by Erd\H{o}s and R\'{e}nyi in~\cite{ErdR} 
the critical probability is equal to $1/n$ (see also~\cite{Bol} or~\cite{JLR} for an 
extensive discussion).

More generally, Bollob\'as, Kohayakawa and \L uczak in~\cite{BKL}
raised the following question: given a sequence of graphs $\{G_n\}$
whose order tends to infinity as $n$ grows, is there such a phase
transition? Assume that $G_n$ has $|G_n|$ vertices and $e_n$ edges.
For each such $n$ we have a probability space on the set of spanning
subgraphs of $G_n$ and the probability of such a subgraph of $G_n$
that has $e$ edges is $p^e(1-p)^{e_n-e}$, where $e_n$ is the number
of edges of $G_n$. Let $G_n(p)$ be a sample from this probability
space. Thus we are seeking a $p_c$ such that: if $p <
(1-\delta)p_c$, then for every $\eps>0$ as $n \rightarrow \infty$
all the components of $G_n(p)$ have at most $\eps |G_n|$ vertices
with probability $1-o(1)$, and if $p > p_c(1+\delta)$, then there
exists $\eps= \eps (p)>0$ for which the largest component of
$G_n(p)$ has at least $\eps |G_n|$ vertices with probability
$1-o(1)$.
If the sequence of graphs is $\{K_n\}$, this is simply the case of a $\np$ random
graph.

Other families of sequences have also been studied.
For example,
percolation on the hypercube with $2^n$ vertices has been analysed
by Ajtai, Koml\'{o}s and Szemer\'{e}di in~\cite{AKS} where
the critical edge probability turns out to be also equal to $1/n$.
More detailed analysis of this phase transition was carried out
recently by  Borgs, Chayes, van der Hofstad, Slade
and Spencer in~\cite{BCHSS2}.

On the other hand recent research has also focused on
finite graphs with bounded maximum degree. Here we consider
sequences of graphs $\{G_n \}$, where each
$G_n$ is a graph on $n$ vertices, with uniformly
bounded maximum degree.
Alon, Benjamini and Stacey have investigated
percolation on such sequences of graphs in~\cite{ABS}. Among other things, they proved that 
the critical probability for the emergence
of a component of linear size in a $d$-regular graph on $n$ vertices
whose girth tends to infinity with $n$
is $1/(d-1)$ (Theorem 3.2 in~\cite{ABS}).
Phase transitions on specific sequences of finite graphs 
were studied more closely by Borgs, Chayes, van der Hofstad, Slade
and Spencer in~\cite{BCHSS},~\cite{BCHSS1}.

More recently, in~\cite{BBCR}, Bollob\'as, Borgs,  Chayes and Riordan analysed the
phase transition in sequences of dense graphs that are convergent in a certain sense.

Also, in ~\cite{FKM} Frieze, Krivelevich and Martin, proved that $p_c=1/d$
for sequences of $d$-regular graphs on $n$ vertices which are quasi-random, when $d \rightarrow \infty$ as $n$ grows.
These are graphs whose structure resembles that of a $d$-regular random graph.

In the present paper, we determine a percolation threshold in the
case where the sequence $\{G_n \}_{n \in \mathbb{Z}^+}$ is a sequence
of sparse random graphs on $n$ vertices. In particular, for every integer $n
\geq 1$, $G_n$ is a uniformly random graph on the set $V_n = \{1,
\ldots, n\}$ having a given degree sequence $\dd (n) = (d_1, \ldots ,
d_n)$, i.e. for $i=1,\ldots n$ vertex $i$ has degree $d_i$.
More formally, a {\em degree sequence} on the set $V_n$ is a vector $\dd =
(d_1,\ldots, d_n)$ consisting of natural numbers, where $d_1 \leq \cdots \leq
d_n$, and $\sum_{i=1}^{n}d_i$ is even.
We let $2M$ denote this sum, and $M=M(n)$ is
the number of edges that $\dd$ spans.
For a given $\dd=\dd(n)$, if $\dd (n)=(d_1,\ldots , d_n)$
for $n \in \mathbb{Z}^+$, we set $D_i = D_i(n)= |\{j\in V_n\ : \ d_j =
i \}|$, for $i\in \mathbb{N}$ and
$\Delta = \Delta(n)=\max_{1\leq i\leq n} \{ d_i\}=d_n$.
Finally, if $G$ is a
graph on $V_n$, then $D(G)$ denotes its degree sequence.

An {\em asymptotic degree sequence} is a sequence  $(\dd (n))_{n \in
\zi^+}$, where for each $n \in \zi^+$ the vector
$\dd (n)$ is a degree sequence on $V_n$. An asymptotic degree
sequence is {\em sparse}, if for every $i \in \mathbb{N}$, we have
$\lim_{n \rightarrow \infty} D_i(n)/n = \la_i$, for some $\la_i \in
[0,1]$, where $\sum_{i\geq 0} \la_i = 1$,  and moreover
\begin{equation} \label{Assumption1}
\lim_{n \rightarrow \infty} {1\over n}\sum_{i\geq 1} i(i-2)D_i(n) = \sum_{i\geq 1}i(i-2)\la_i < \infty.
\end{equation}
This implies that for every $\eps > 0$ there exists $i^*(\eps)$ and $N = N (\eps)$ such that
for every $n > N$ we have
\begin{equation} \label{Assumption2}
\left|{1\over n}\sum_{i\leq i^*} i(i-2)D_i(n) - \sum_{i\geq 1}i(i-2)\la_i  \right| < \eps.
\end{equation}
The {\em generating polynomial} of a sparse asymptotic degree
sequence is defined as $L(s)=\sum_{i=0}^{\infty}\la_i s^i$. We assume that
every asymptotic degree sequence $(\dd (n))_{n \in \mathbb{Z}^+}$ we
work with is such that for every $n$ the set of simple graphs that
have $\dd(n)$ as their degree sequence is non-empty.

We consider two types of percolation. Firstly, for some $p \in
(0,1)$, each edge of $G_n$ is present with probability $p$
independently of every other edge. This type of percolation is usually
called {\em bond percolation}, in that we randomly delete the edges
(i.e. the bonds) of $G_n$. This is distinguished from another type of
percolation which is called {\em site percolation}. Here, we go
through the vertices of $G_n$ and we make each of them isolated with
probability $1-p$, independently of every other vertex (or as we say we
delete this vertex). The random subgraph in this case is the spanning
subgraph of $G_n$ that does not contain the edges that are attached to
the vertices that were deleted. The
terms  ``bond'' and ``site'' percolation have their origins in the
percolation theory of infinite graphs
(see~\cite{Grim} for an extensive discussion on both types as well
as the references therein).

We shall now define the percolation threshold in each of the above
cases. Let $G'(n)$ denote the random subgraph that is obtained in
either case and let $L_1 (G'(n))$ be the lexicographically first
component of $G'(n)$ (this is the component that has maximum order
and the smallest vertex it contains is smaller than the smallest
vertex of every other component of maximum order - the comparison
between the vertices is by means of the total ordering on $V_n$).
Starting from the bond percolation we set
$p_c^{bond} = \sup \{ p \in [0,1] : \ |L_1 (G'(n))|/n \stackrel{p}
\rightarrow  0  \ \mbox{as $n \rightarrow \infty$} \}$ (the symbol $\stackrel{p}
\rightarrow$ denotes convergence in probability, i.e. we say that $X_n \stackrel{p}
\rightarrow 0$ if for every $\eps>0$ we have
$\p [|X_n|> \eps] \rightarrow 0 $ as $n\rightarrow 0$). The
convergence in probability is meant with respect to the sequence of probability
spaces indexed by the set $\mathbb{Z}^+$, where for each
$n \in \mathbb{Z}^+$ the probability of a certain spanning subgraph is
the probability that this is the subgraph which is spanned by the
edges that survive the random deletion of the edges of the random graph
$G_n$. Similarly, in the case of site percolation we define $p_c^{site} =
\sup \{ p \in [0,1] : \ |L_1 (G'(n))|/n \stackrel{p}
\rightarrow  0  \ \mbox{as $n \rightarrow \infty$} \}$, where $G'(n)$
is now the spanning subgraph of $G_n$ that is the outcome of the deletion of
those edges that attached to the chosen vertices, i.e. the vertices
that we make isolated.
Note that in both cases
there are two levels of randomness.

If $G_n$ is a random $d$-regular graph on $V_n$, for
any fixed $d\geq 3$, the bond percolation threshold has been
calculated by Goerdt in~\cite{Goe} and is equal to $1/(d-1)$.
Before this, bond percolation in random regular graphs was studied
by Nikoletseas, Palem, Spirakis and Yung
in~\cite{NikS1}, where it was proved that the critical probability
is at most $32/d$, for $d$ large enough. Also, Nikoletseas and Spirakis in~\cite{NikS}
study the edge expansion properties of the giant component that remains after the
edge deletion process.
However, these papers did not provide any analysis
on the site percolation process.
Our main theorem involves also the latter and is stated as follows:
\begin{theorem} \label{main}
If $(\dd (n))_{n \in \mathbb{Z}^+}$ is a sparse asymptotic degree sequence
of maximum degree $\Delta (n) \leq n^{1/9}$
and $L(s)$ is its generating polynomial which is twice differentiable
at 1 and moreover $L''(1) > L'(1)$, then
$p_c^{site}=p_c^{bond}=L'(1)/L''(1)$. Moreover, whenever
$p>p_c^{bond}$ ($p>p_c^{site}$, respectively) there exists an
$\varepsilon >0$ such that
$|L_1(G'(n))|>\varepsilon n$ with probability $1-o(1)$.
\end{theorem}
The formula for both critical probabilities was obtained by Dorogovtsev
and Mendes in~\cite{DorMend} using qualitative (i.e. non-rigorous) arguments.

To make the statement of the above theorem slightly clearer, let us
consider the case of bond percolation  (the case of site
percolation is similar). Let $\mathcal{G}(n)$ be the set of graphs
on $V_n$ whose degree sequence is $\dd (n)$. Each graph $G \in
\mathcal{G}(n)$ gives rise to a probability space which consists of
all its spanning subgraphs. In particular, if $G$ has $e$ edges and
$G'$ is a spanning subgraph of $G$ that has $e' \leq e$ edges then
its probability is $p^{e'}(1-p)^{e-e'}$; let $\p_p^G[\cdot ]$ denote
this measure. In other words, this space accommodates the outcomes
of the bond percolation process applied to $G$ and we call it the
\emph{percolation space} of $G$. For any $\varepsilon \in (0,1)$ we
let $g_{\varepsilon}(G)$ be the set of all spanning subgraphs of $G$
whose largest component has at least $\varepsilon n$ vertices. This
event has probability $\p_p^G[g_{\varepsilon}(G)]$ in the
percolation space of $G$.
Now, assume that $p < p_c^{bond}$. Theorem~\ref{main} implies that for any given $\rho \in
(0,1)$, the event
$\{G \in \mathcal{G}(n)\ :\  \p_p^{G}(g_{\varepsilon}(G))< \rho\}$
occurs with probability
$1-o(1)$ in the uniform space $\mathcal{G}(n)$.
That is, asymptotically for almost every graph in
$\mathcal{G}(n)$ the random deletion of the edges leaves a
component of order at least $\eps n$
with probability no more that $\rho$.
If $p>p_c^{bond}$, then the second part of the theorem implies that
there exists $\varepsilon >0$ such that the event $\{G \in
\mathcal{G}(n)\ : \ \p_p^{G}(g_{\varepsilon}(G))>1-\rho \}$ occurs
with probability $1-o(1)$ in $\mathcal{G}(n)$. Hence, as
$n\rightarrow \infty$ almost all the graphs in $\mathcal{G}(n)$ are
such that if we apply the bond percolation process to them with
retainment probability $p$, then there is a component having at
least $\eps n$ vertices with probability at least $1-\rho$ (in the
percolation space).

The fact that the critical probabilities coincide reflects a behaviour
that is similar to that of percolation on an infinite regular
tree. Of course in that context
the critical probabilities are defined with respect to the appearance
of an infinite component that contains the (vertex that has been selected as the) root.
Using the fundamental theorem of Galton-Watson processes (see for example~\cite{Branch}),
it can be shown that the bond and the site
critical probabilities coincide and they are equal to $1/(d-1)$, where
$d$ is the degree of each vertex of the tree. Observe that for the case
of a random $d$-regular graph the above theorem implies that
$p_c^{site}=p_c^{bond}=1/(d-1)$.
This is not a coincidence as it is
well-known that a random $d$-regular graph  locally
(e.g. at distance no more than $i$ from a given vertex for some fixed $i$)
looks like  a $d$-regular tree.

More generally, the typical local structure of the class of
random graphs we are investigating are also tree-like. 
Note that the ratio $L''(1)/L'(1)$ equals
\begin{equation} \label{NoOfChildren}
\sum_{i=2}^{\infty}(i-1)\frac{i\la_i}{\sum_{j=1}^{\infty}j\la_j}.
\end{equation}
Consider a vertex $v\in V_n$ which has positive degree and let us
examine more closely the behaviour of one of its neighbours. It can
be shown that the probability that this has degree $i$ is
proportional to $i D_i(n)$. In particular, it is almost equal to
$\frac{iD_i(n)}{\sum_i i D_i(n)}$ and this tends to
$\frac{i\la_i}{\sum_{j=1}^{\infty}j\la_j}$ as $n$ grows. Moreover,
one can show that with probability $1-o(1)$ there are no edges
between the neighbours of $v$. Therefore (\ref{NoOfChildren}) is the
limit of the expected number of children a neighbour of $v$ has.
This scenario is repeated for every vertex in the $d$-th
neighbourhood of $v$, where $d$ is fixed. More precisely, the
vertices which are at distance no more than $d$ induce a tree rooted
at $v$ which contains at most $\ln\ln~n$ vertices, with probability
$1-o(1)$. Suppose that there are $t_i$ vertices of degree $i$ in
this tree. Thus for a vertex that is at distance $d$ from $v$, the
probability that it has degree $i$ is proportional to $i
(D_i(n)-t_i)=i\la_i n(1-o(1))$. More precisely, it is
$$ {i (D_i(n)-t_i) \over \sum_i i(D_i(n)-t_i)}.$$
Since $\Delta \leq n^{1/9}$ and $t_i\leq \ln \ln n$, it follows that
$\sum_i it_i \leq \ln \ln n \sum_{i\leq n^{1/9}} i = O(n^{1/3})$.
Hence, the limit of the above probability as $n \rightarrow \infty$ is again
$\frac{i\la_i}{\sum_{j=1}^{\infty}j\la_j}$ and (\ref{NoOfChildren}) gives the
limiting expected number of children of such a vertex.
In other words, the graph that is induced by the
vertices which are at distance no more than $d$ from $v$ behaves like
the tree of a branching process that started at $v$,
with the ratio $L''(1)/L'(1)$ being the expected progeny of each vertex.
Observe here that the condition $L''(1) > L'(1)$ implies that in fact
this is a supercritical branching process which yields an infinite tree
with probability 1.

Therefore, at least locally either bond or site percolation is
essentially percolation on such a random rooted tree.
In both types of percolation, if $p< L'(1)/L''(1)$,
then the expected number of children of a vertex that survive is $p L''(1)/L'(1)<1$.
Thus the random tree that is developing around $v$ after the random failures of the
edges or the vertices
will be distributed as the tree of a subcritical branching process.
In particular, the tree that surrounds most of the
vertices will be cut off from the rest of the graph at a relatively small depth.
On the other hand, if $p> L'(1)/L''(1)$ a large
proportion from each of these local trees is preserved and moreover they are big enough to guarantee that
there are enough edges going out of them. So eventually there is a fair chance that
some of them are joined together and form a component of linear order. However, this is only a
qualitative approach to Theorem \ref{main}. The actual proof
and the structure of the paper are described in Section~\ref{Sketch}.

\subsection{Theorem~\ref{main} and power-law graphs}
The power-law degree sequences are those for which for any $k\geq 1$
one has
$\la_k =c{k^{-\gamma}}$, for some constants $c, \gamma >0$.
We should point out that the crucial parameter here is $\gamma$.
Such degree sequences
have attracted much attention in the last few years mainly because of the fact that they arise
in ``natural'' networks such as the Internet, the WWW or biological networks
(see~\cite{BolR},~\cite{BA} or~\cite{DorMend} for a survey of results or the recent
book by Chung and Lu~\cite{Chung} for a more detailed discussion).
For example, in~\cite{FFF} Faloutsos, Faloutsos and Faloutsos gave evidence that the Internet as it looked like
in 1995, viewed as
a graph whose vertices are the routers, and
the edges are the physical links between them, has a power-law degree sequence with $\gamma \approx 2.48$.
Bond and site percolation on such networks naturally correspond to random failures of the links
or of the nodes, respectively. Thus, a site percolation process on the Internet
may be seen as random failures of routers.

For a power-law degree sequence with $\gamma > 3$ one has 
$L(s) = c\sum_{k\geq 2} {s^k \over k^\gamma}$
where $|s|\leq 1$. Thus, if $\zeta (\lambda )= \sum_{k \geq 1}{1\over k^\lambda}$ is the Riemann's zeta
function, then 
$L'(1) = c\sum_{k \geq 2} {k \over k^\gamma} =c\sum_{k \geq 2} {1 \over k^{\gamma-1}} =c\zeta (\gamma -1)$ and 
$L''(1) = c\sum_{k \geq 2} {k (k-1)\over k^\gamma} = c\sum_{k \geq 2} {1\over k^{\gamma-2}} - c\sum_{k \geq 2} {1\over k^{\gamma-1}}
= c(\zeta (\gamma -2)-\zeta (\gamma -1))$ (of course here the derivatives are left derivatives).
Let $\gamma_0 = \sup \{ \gamma\ : \ \gamma >3, \ {\zeta(\gamma-2) \over \zeta(\gamma-1)}> 2 \}$.
Theorem~\ref{main} implies that the critical probabilities for a power-law degree sequence with
$3 < \gamma <\gamma_0$ but maximum degree at most $n^{1/9}$ are
\begin{equation}\label{PowerLaw}
p_c^{site} = p_c^{bond} = {\zeta (\gamma -1)\over \zeta(\gamma-2) - \zeta(\gamma-1)}.
\end{equation} 
This agrees with the analysis made in~\cite{CNSW} by
Callaway, Newman, Strogatz and Watts for the case of site
percolation on a random graph whose degree sequence follows a "truncated"
power-law, that is $\la_k = Ck^{-\gamma} e^{-k/ \kappa}$, for $C,
\kappa >0$. As $\kappa \rightarrow \infty$, then this approaches a
power-law distribution with parameter $\gamma$. It can be shown that
in this case the critical probability they obtain converges to the
expression in (\ref{PowerLaw}) (see for example Equation (141) p.45
in~\cite{BA}). Similar analysis made by Cohen, Erez, ben-Avraham and
Havlin in~\cite{CEbAH} suggests that if $\gamma \leq 3$ there is no
phase transition at all. Also, Bogu\~n\'a, Pastor-Satorras and
Vespignani in~\cite{BoPaVe} argue that this happens whenever $2<
\gamma \leq 3$. This was also suggested by simulations in ~\cite{AJB}. 
In particular, Albert, Jeong and Barab\' asi give experimental
evidence of the result of a site percolation process on a random
graph (obtained from a different model) whose degree sequence is power-law with $\gamma =3$. They
observe that the graph remains largely intact by the random deletion
of vertices and no threshold behaviour is observed. That is, even if
a large proportion of vertices are deleted, there is always a
component of linear order. In the same paper, the authors give also
experimental evidence in samples of the Internet and the World-Wide
Web, concluding that no phase transition occurs even for small
values of $p$. However, Dorogovtsev and Mendes 
in~\cite{DorMend} applying the formula of Theorem~\ref{main} (which they also obtain in their paper), 
but without our degree restrictions, 
give the scaling of the critical probabilities as functions of $n$ as $n\rightarrow \infty$, 
for $2<\gamma \leq 3$. In that context the critical probabilities are defined empirically,  
according to whether or not the proportion of vertices in the largest component after the percolation process is 
almost zero.

The case $\gamma \leq 3$ corresponds to $L''(1)$ being divergent which suggests that
$p_c^{site}$ and $p_c^{bond}$ vanish.
However, Theorem~\ref{main}
works under the assumption that $L''(1)$ converges.
It would be an interesting and natural next step to prove or disprove the existence of a
positive critical probability in the case where $L''(1)$ is divergent.

\section{Definitions and sketch of the proof} \label{Sketch}
In this paper we are interested in sparse asymptotic degree sequences
$\mathcal{D}$
satisfying
$Q(\mathcal{D}):=\sum_{i=1}^{\infty} i(i-2)\la_i >0$. This is
equivalent to saying that $\sum_{i=1}^{\infty} i(i-1)\la_i =
L''(1)>L'(1) = \sum_{i=1}^{\infty} i \la_i$, where $L(s)$ is the
generating function of $\mathcal{D}$.

One of the main tools we use in the present work is the configuration
model, which was introduced in different versions 
by Bender and Canfield in $\cite{BCan}$ and Bollob\'as in~\cite{Bol1}.
If $\dd$ is a degree sequence on $V_n$, for some $n \in \zi^+$,
we define the set of points $P=P(\dd)$ as $\{1\times [d_1], \ldots, n \times
[d_n]\}$, where $[d_i]=\{1,\ldots, d_i \}$ if $d_i >0$ or the empty
set otherwise. That is to every vertex in $V_n$ correspond $d_i$
points. Clearly, there are $2M$ points in $P$.
Thus observe that there are ${(2M)! \over M!2^{M}}$ perfect matchings on $P$.
If $M(\dd)$ is such a perfect matching, then we can obtain a (multi)graph $GM(\dd)$
if we project $P$ onto $V_n$ preserving adjacencies, namely for any
two vertices $i,j \in V_n$, if $M(\dd)$ contains an edge between a
point in $i \times [d_i]$ and a point in $j \times [d_j]$, then
$GM(\dd)$ contains a copy of the edge $(i,j)$. Of course in a perfect
matching there might be edges that join two points corresponding to the same
vertex, in which case $GM(\dd)$ obtains a loop on the
vertex. Similarly, there might be two vertices which are joined to
each other with more than one pairs of points, and in this case $GM(\dd)$
obtains multiple copies of the corresponding edge. If $M(\dd)$ is a
uniformly random perfect matching on $P$, then observe that $GM(\dd)$
is not uniformly distributed over the set of multigraphs having $\dd$
as their degree sequence. However, if we condition on the event that
$GM(\dd)$ is a simple graph, then it is uniformly distributed over the
set of simple graphs that have $\dd$ as their degree sequence.

Consider now an asymptotic degree sequence
$\mathcal{D}=(\dd (n))_{n \in \zi^+}$. For each $n \in \zi^+$ we set
$P(n)$ to be the set of points that corresponds to the degree sequence
$\dd (n)$. Let $M_n$ be a uniformly random perfect matching on $P(n)$
and let $\tilde{G}(n)$ be the multigraph that is obtained from the projection of
$M_n$ onto $V_n$. The following theorem was proved by M. Molloy and
B. Reed in \cite{MoR} and has a key role in our proofs:
\begin{theorem} \label{key}
Let $\mathcal{D}=(\dd (n))_{n \in \zi^+}$ be a sparse
asymptotic degree sequence of maximum degree at most $n^{1/9}$.
\begin{itemize}
\item If $Q(\mathcal{D})>0$, then there exists an $\varepsilon >0$
such that $\p [|L_1(\tilde{G}(n))|\geq \varepsilon n]\rightarrow 1$, as $n
\rightarrow \infty$.
\item If  $Q(\mathcal{D})<0$, then for every $\varepsilon >0$ we have
$\p [|L_1(\tilde{G}(n))|\geq \varepsilon n]\rightarrow 0$, as $n
\rightarrow \infty$.
\end{itemize}
\end{theorem}
Of course the above theorem as stated in~\cite{MoR} was referring to
simple graphs rather than multigraphs. However, it was actually
proved for the random multigraph $\tilde{G}(n)$ and conditioning on
being simple it can be stated for random simple graphs having this
particular degree sequence. In fact in the first case Molloy and
Reed proved the uniqueness of the component that has linear order;
in particular the second largest component has logarithmic order.
The restriction on the maximum degree can be slightly relaxed
(see~\cite{MoR}), but for the simplicity of our proofs we assume it
to be as above. The way we use this result will become apparent
during the sketch of our proofs that is about to follow.


Note that assuming that $\mathcal{D}$ has $L''(1)>L'(1)$, the
above theorem implies that $\tilde{G}(n)$ will have a giant component with
probability $1-o(1)$.

Here are the two deletion processes that we consider separately:
\begin{itemize}
\item {\em bond percolation:} For some $p \in (0,1)$, we delete at
random each edge of $\tilde{G}(n)$ with probability $1-p$, independently of
every other edge.
\item {\em site percolation:} with probability $1-p$ we make a vertex
isolated by deleting the edges that are incident to it, independently
for every vertex of $\tilde{G}(n)$.
\end{itemize}
In either case,
the random multigraph that is the outcome of this experiment is denoted by
$G'(n)$.

Eventually we want to know the structure of $G'(n)$, if $\tilde{G}(n)$ is a simple graph.
Thus, we will show that
if $\mathcal{A}(n)$ is set of multigraphs on $V_n$
and $\p [G'(n) \in \mathcal{A}(n)]\rightarrow 0$
as $n \rightarrow \infty$, then $\p [G'(n) \in
\mathcal{A}(n) \ | \ \tilde{G} (n) \mbox{ is simple}]\rightarrow 0$ as $n \rightarrow \infty$ as
well. Hence, it will be sufficient for our purposes to perform the
random deletion on the edges or the vertices of a
random perfect matching on $P(n)$ without any conditioning and
henceforth to consider the multigraph that is obtained from the
remaining edges; this is going to make the calculations much simpler.

Thus we prove the following lemma:
\begin{lemma} \label{simple}
Let $\mathcal{A}(n)$ be a set of multigraphs on $V_n$ and suppose that
$\p [G'(n) \in \mathcal{A}(n)]\rightarrow 0$ as $n \rightarrow \infty$.
Then $\lim_{n \rightarrow \infty} \p [G'(n) \in \mathcal{A}(n) \ | \ \tilde{G} (n) \mbox{ is simple}]= 0$ as
well.
\end{lemma}
\begin{proof}
Note that
\begin{eqnarray} \label{SimpleGraphs}
\p [ G'(n) \in \mathcal{A}(n)\ |
\ \mbox{$\tilde{G}(n)$ is simple}]  &\leq & \frac{\p [ G'(n) \in \mathcal{A}(n)]}
{\p [\mbox{$\tilde{G}(n)$ is simple}]}.
\end{eqnarray}
The asymptotic enumeration formula for graphs with a given degree sequence, such 
that $M=\Theta(n)$   
and $\Delta = o(n^{1/3})$, obtained by McKay and Wormald in~\cite{McW}
yields
$$\p[ \mbox{$\tilde{G}(n)$ is simple}] =
(1+o(1)) e^{-\la/2-\la^2/4},$$
where $\la = {1\over M} \sum_{i=1}^{n} {d_i \choose 2}$.
But the latter sum is at most $\sum_{i=1}^{n} d_i^2 = O(n)$, and
since $M=\Omega (n)$ it follows that $\la = O(1)$.
Thus
$$ \liminf_{n\rightarrow \infty} \p[ \mbox{$\tilde{G}(n)$ is simple}]  >0,$$
and this concludes the proof of the lemma as the numerator in
(\ref{SimpleGraphs}) converges to zero.
\end{proof}

In both cases the random deletion induces a (random) degree sequence
on $V_n$ which we denote by $\dd'(n)$ (the use the same symbol
for the two kinds of percolation should cause no confusion). So, for each $n
\in \zi^+$, let $\D_n$ be the set of degree sequences on $V_n$ that
are the result of the random deletion equipped with the probability
distribution inherited by the random experiment we just described.
That is, the probability of a certain degree sequence $\dd'(n)\in \D_n$ is
the probability that the degree sequence which is induced by the random
deletion (either of edges or of vertices) on $\tilde{G}(n)$ is $\dd'(n)$.
We set $\D = \prod_{n=1}^{\infty}\D_n$ to be the product space
equipped with the product measure, which
we denote by $\mu$. Thus each element of $\D$ is an asymptotic degree
sequence and $Q$ is now a random variable on $\D$.

The strategy of our proof is quite different from that in~\cite{Goe}, in that we make
explicit use of Theorem~\ref{key}.
We first prove that the perfect matching between those points that are the endpoints
of the edges that survive the deletion either in bond percolation or
in site percolation is uniformly random among the perfect matchings on
these points.
Hence to study the asymptotic properties
of $G'(n)$ we shall condition first on its degree sequence for every $n
\in \zi^+$ and then we shall study the asymptotic behaviour of $G'(n)$
conditioned on this asymptotic degree sequence. Of course to show that $G'(n)$
has a certain property with probability tending to 1 as $n \rightarrow
\infty$, we have to show that almost all the asymptotic degree
sequences in $\D$ have similar behaviour. In particular, if $D_i'(n)$ is the number of 
vertices of degree $i$ in $G'(n)$, we shall 
prove that the random variable ${1\over n} \sum_{i\geq 1}i(i-2)D_i'(n)$ converges $\mu$-almost surely ($\mu$-a.s.)
to a quantity $Q'$ that depends only on the $\la_i$'s and on $p$,
which we will calculate explicitly in both cases.
From this we derive the critical $p_c$, which we denote by $p_c^{bond}$
for the case of bond percolation  and $p_c^{site}$ for the case of
site percolation. We show that if $p>p_c$ then $Q'$ is positive, 
whereas if $p<p_c$ we have $Q'< 0$. Using Theorem \ref{key}, we deduce the sudden
appearance of a giant component in $G'(n)$ when $p$ crosses $p_c$,
with probability that tends to 1 as $n \rightarrow \infty$.

We conclude this section stating a concentration
inequality which we use in our proofs and it follows from Theorem 7.1 in~\cite{McD}.
Let $S$ be a finite set and let $f$ be a real-valued
function on the set of those subsets of $S$ that have size $k$. Assume
that whenever $c, c'$ are two such subsets whose symmetric difference is
$2$, then $|f(c)-f(c')|\leq 2$. If $C$ is chosen uniformly at
random among the $k$-subsets of $S$, then
\begin{equation} \label{Conc}
\p \left[|f(C)- \ex[f(C)]|>t\right] \leq 2 \exp \left(-{t^2 \over 2k} \right).
\end{equation}

\section{Bond percolation}
In this case, we start with the random graph $\tilde{G}(n)$ which the multigraph that is the
projection onto $V_n$ of a uniformly random perfect matching on
$P(n)$ and we create the multigraph $G'(n)$, deleting each edge of the matching with probability
$1-p$, independently of every other edge. Thus, the number of edges
of $G'(n)$ is distributed as $\mathrm{Bin}(M(n),p)$.

Firstly, we will prove that the perfect matching on the
remaining points in $P(n)$ conditional on
the degree sequence that is created after the deletion
is uniformly distributed on the set of perfect
matchings on the set of points in $P(n)$ that survive the deletion.
In particular, if $C$ is the set points in $P(n)$ that are the
end-points of the surviving edges, then for every $i \in V_n$
the new degree of vertex $i$ is $|C\cap (i\times [\dd_i(n)])|$.
Hence the (random) set $C$ induces a degree sequence on $V_n$, which we
denote $\dd'(n)$. We set $P'(n)=P(\dd'(n))$.

Let $d \in \{ 0,\ldots , \Delta \}$ and assume that the vertices
$i_1, \ldots , i_{k_d}$ (and only these) have new degree
equal to $d$ after the edge deletion.
Hence $\dd'(n)$ contains exactly
$k_d$ vertices of degree $d$ and assume that these are $i,\ldots, i+k_d-1$.
We identify $i_j$ with $i+j-1$, for every $j=1, \ldots,
k_d$. Moreover, provided that $d\geq 1$ we also identify the $d$ points of
$C\cap (i_j\times [d(i_j)])$ with the points $\{i+j-1\} \times \{1,\ldots ,d
\}$ in $P'(n)$.
Hence any perfect matching between the points in $C$
corresponds to a perfect matching on $P'(n)$ and vice-versa. In other words,
we obtain a bijection between the perfect matchings on these two
sets of points.

In the case of bond percolation, the set $\D_n$
consists precisely of those degree sequences that are induced by the
deletion of the edges of a random perfect matching on $P(n)$.
The probability of a certain degree sequence in
$\D_n$ will therefore be the probability that this is the induced
degree sequence after the deletion.
Our aim is to show that conditional on $\dd'(n)=\dd'$, each perfect matching on
$P'(n) = P(\dd')$ is equilikely.

To do so, we first prove the following:
\begin{lemma} \label{uni-surv}
Conditional on having $k$ edges that survive the random deletion of the edges
of the perfect matching on $P(n)$, the set of their $2k$ end-points is
uniformly distributed among the $2k$-subsets of $P(n)$.
\end{lemma}
\begin{proof}
The probability that a specific $2k$-subset of $P(n)$ is the set of
the end-points of the $k$ edges that survive is the probability that
the perfect matching on $P(n)$ consists of a perfect matching on these $2k$
points and a perfect matching between the $2M-2k$ remaining points,
and that it is the set of the $k$ edges on this $2k$-subset that
survive the deletion. The probability of this event is exactly:
\begin{eqnarray*}
\frac{\frac{(2M-2k)!}{(M-k)!2^{M-k}}\frac{2k!}{k!2^k}}{\frac{2M!}{M!2^M}}
\ \frac{1}{{M \choose k}} = \frac{1}{{2M \choose 2k}}.
\end{eqnarray*}
This concludes the proof of the lemma.
\end{proof}
With a little more work
we obtain what we were aiming for:
\begin{lemma} \label{uniform}
Let $\dd'(n)$ be the degree sequence that is induced by the random
deletion of the edges of a uniformly random perfect matching on
$P(n)$. Conditional on $\dd'(n) =\dd'$, any perfect matching on
$P(\dd')$ is equilikely.
\end{lemma}
\begin{proof}  Assume that the sum of the degrees in $\dd'$ is $2k$ and
let $S_{\dd'}$ be the set of those $2k$-subsets of $P(n)$ that
induce the degree sequence $\dd'$.
Let $m$ be a particular perfect matching on $P'(n)$, conditional on $\dd'(n)=\dd'$.
In other words, $m$ is a perfect matching on $P(\dd')$.
Now, let us condition on $|P'(n)|=2k$.
If $C' \in S_{\dd'}$, the probability that $C=C'$
and the particular perfect matching that corresponds to $m$ is realised on $C$ is
${1\over {M \choose k}}~ \frac{\frac{(2M-2k)!}{(M-k)!2^{M-k}}}{\frac{2M!}{M!2^M}}$.
Thus,
\begin{equation} \label{1stcond}
\p [m, \ \dd'(n)=\dd' \ | \ |P'(n)|=2k] = \sum_{C' \in S_{\dd'}}
{1\over {M \choose k}}~ \frac{\frac{(2M-2k)!}{(M-k)!2^{M-k}}}{\frac{2M!}{M!2^M}}
= {|S_{\dd'}|\over {M \choose k}}~ \frac{\frac{(2M-2k)!}{(M-k)!2^{M-k}}}{\frac{2M!}{M!2^M}}.
\end{equation}
By the previous lemma, conditional on $|P'(n)|=2k$ every set in $S_{\dd'}$ has probability
$1/{2M \choose 2k}$. Therefore,
$$ \p [\dd'(n)=\dd' \ | |P'(n)|=2k ] = \frac{|S_{\dd'}|}{{2M \choose 2k}}.$$
Now, Bayes' rule (i.e. dividing (\ref{1stcond}) by the above probability) yields:
\begin{eqnarray*}  \p [m \ | \ \dd'(n)=\dd', \ |P'(n)|=2k] &=&
{1\over {M \choose k}}~ \frac{\frac{(2M-2k)!}{(M-k)!2^{M-k}}}{\frac{2M!}{M!2^M}} ~ {2M \choose 2k}\\
&=& \frac{k!(M-k)!}{M!} \
\frac{\frac{(2M-2k)!}{(M-k)!2^{M-k}}}{\frac{2M!}{M!2^M}} \
\frac{2M!}{(2M-2k)!2k!}= \frac{1}{\frac{2k!}{k!2^k}}.
\end{eqnarray*}
But
\begin{eqnarray*}
\p [m \ | \ \dd'(n)=\dd'] &= & \p [m \cap |P'(n)|=2k \ | \ \dd'(n) = \dd']  \\
&=&
\p [m \ | \  |P'(n)|=2k, \ \dd'(n) = \dd']~ \p [|P'(n)|=2k \ | \ \dd'(n) = \dd'] \\
&=& \p [m \ | \  |P'(n)|=2k, \ \dd'(n) = \dd'],
\end{eqnarray*}
since $\p [|P'(n)|=2k \ | \ \dd'(n) = \dd'] = 1$, and the lemma follows.
\end{proof}


For $i\in \mathbb{N}$, let $D_i'(n)$ be the number of vertices of
degree $i$ in $\dd'(n)$.
The threshold probability will be determined by the quantity
$\sum_{i\geq 1} i(i-2)\la_i^{bond}$, where
$$\la_i^{bond} = \lim_{n \rightarrow \infty} {1\over n} \ex[D_i'(n)].$$
Hence we need to determine each $\la_i^{bond}$, proving the existence
of this limit, and to do so
we will first calculate the expected value of
$D_i'(n)$. We begin with the conditional
expectation given the size of $C$, that is conditional on the number
of points in $P(n)$ that survive the deletion. For any $k=0,\ldots,
M$ we have
\begin{eqnarray*}
\ex [D_i'(n) \ | \ |C|=2k] = \sum_{d=i}^{\Delta} D_d (n)
\p [\mbox{a given vertex of degree $d$ has new degree $i$}\ | \ |C|=2k].
\end{eqnarray*}

Recall that $|C|/2$, which equals the number of edges that survive
the random deletion,
is distributed as
$\mathrm{Bin}(M,p)$. Therefore, a standard concentration argument
yields
\begin{equation} \label{CConc}
\p [||C|/2 - Mp|> \ln n \sqrt{n}]\leq \exp
\left(-\Omega( \ln^2 n )\right).
\end{equation}
This indicates that we may
restrict ourselves to $k \in I=[Mp -\ln n\sqrt{n}, Mp + \ln n \sqrt{n}]$.

By Lemma~\ref{uni-surv} conditional on $|C|=2k$, the set $C$ is
uniformly distributed among all $2k$-subsets of $P(n)$.
Hence, we obtain
\begin{eqnarray*}
\lefteqn{\p [\mbox{a given vertex of degree $d$ has new degree
$i$}\ | \ |C|=2k]=}\\
& & {d \choose i} \frac{{2M-d \choose 2k-i}}{{2M \choose 2k}} =
{d \choose i} \frac{(2M-d)!}{2M!} \ \frac{2k!}{(2k-i)!} \
\frac{(2M-2k)!}{(2M-d-2k+i)!} \\
&=& {d \choose i} \frac{(2k)^i}{(2M)^d} \ (2M-2k)^{d-i}
\left(1+ O\left(\frac{1}{n^{7/9}} \right) \right) =
{d \choose i} p^i (1-p)^{d-i} \left(1 + O\left(\frac{\ln n}{n^{7/18}} \right) \right),
\end{eqnarray*}
uniformly for any $d \leq \Delta$ and any $k\in I$.

Therefore, since $D_d'(n) \leq n$ by (\ref{CConc}) we obtain:
\begin{eqnarray*}
\ex[D_i'(n)] &=& \sum_{k=0}^{M}\ex [D_i'(n) \ | \ |C|=2k] \p[|C|=2k] \\
&=&
\sum_{k\in I}
\sum_{d=i}^{\Delta} D_d (n){d \choose i} p^i
(1-p)^{d-i}  \left(1+ O\left(\frac{\ln n}{n^{7/18}} \right) \right) \p[|C|=2k]
\\& &  + o\left(\frac{1}{n^3} \right) \\
&=& \left(1+  O\left(\frac{\ln n}{n^{7/18}} \right) \right)
\sum_{d=i}^{\Delta} D_d (n){d \choose i} p^i
(1-p)^{d-i}
+ o\left(\frac{1}{n^3} \right).
\end{eqnarray*}
For every $\eps >0$, if $i'$ and $n$ are large enough
$${1\over n} \sum_{d=i'+1}^{\Delta} D_d (n){d \choose i} p^i(1-p)^{d-i} \leq
{1\over n} \sum_{d=i'+1}^{\Delta} D_d (n) < \eps .  $$
Therefore,
$${1\over n} \sum_{d=i}^{i'} D_d (n){d \choose i} p^i
(1-p)^{d-i} \leq {1\over n}\sum_{d=i}^{\Delta} D_d (n){d \choose i} p^i
(1-p)^{d-i} \leq {1\over n} \sum_{d=i}^{i'} D_d (n){d \choose i} p^i
(1-p)^{d-i} + \eps. $$
Taking limits on both sides we obtain:
$$  \sum_{d=i}^{i'} \la_d {d \choose i} p^i
(1-p)^{d-i} \leq \liminf_{n \rightarrow \infty}{1\over
n}\sum_{d=i}^{\Delta} D_d (n){d \choose i} p^i (1-p)^{d-i} $$
and
$$\limsup_{n \rightarrow \infty} {1\over n}\sum_{d=i}^{\Delta} D_d (n){d \choose i} p^i
(1-p)^{d-i} \leq \sum_{d=i}^{i'} \la_d {d \choose i} p^i
(1-p)^{d-i} + \eps.  $$
Letting $i'\rightarrow \infty$ and then $\eps \rightarrow 0$, we obtain
the value of $\la_i^{bond}$:
\begin{equation} \label{labond}
\lim_{n \rightarrow \infty} \frac{\ex[D_i'(n)]}{n} =
\sum_{d=i}^{\infty} \la_d {d \choose i} p^i (1-p)^{d-i} \equiv \la_i^{bond}.
\end{equation}
We will show that
the critical probability $p^{bond}_c$ is equal to the root
of the equation $Q':=\sum_{i=1}^{\infty}i(i-2)\la_i^{bond}=0$, which
we denote by $\hat{p}_{bond}$. Firstly, let us calculate
$\hat{p}_{bond}$.
We have
\begin{eqnarray*} \label{threshold}
\sum_{i=1}^{\infty}i(i-2)\la_i^{bond}&=& \sum_{i=1}^{\infty}
i(i-2)\sum_{d=i}^{\infty} \la_d {d \choose i}p^i (1-p)^{d-i} =
\sum_{d=1}^{\infty} \la_d \sum_{i=1}^{d} i(i-2) {d \choose i}p^i
(1-p)^{d-i} \nonumber \\
&=& \sum_{d=1}^{\infty} \la_d
\left(dp(1-p)+(dp)^2-2dp \right) = \sum_{d=1}^{\infty} \la_d
\left(1-p +dp-2 \right)dp \nonumber \\
& =& \sum_{d=1}^{\infty} \la_d
\left((d-1)p-1 \right)dp
= \sum_{d=1}^{\infty} \la_d d(d-1)p^2 - p \sum_{d=1}^{\infty} d
\la_d \nonumber \\
&=& p \left(p \sum_{d=1}^{\infty} \la_d d(d-1) -  \sum_{d=1}^{\infty} d
\la_d  \right).
\end{eqnarray*}
Therefore,
$$\hat{p}_{bond} := \frac{L'(1)}{L''(1)}.$$

We now let $Q_n' = {1\over n}\sum_{i\geq 2} i(i-2)D_i'(n)$ and will
show that $\lim_{n \rightarrow \infty} Q_n'$ exists $\mu-a.s.$ and it
is equal to $Q'$. Hence, the sign of $Q'$ determines the sign of
$\lim_{n \rightarrow \infty} Q_n'$ for almost every asymptotic degree
sequence in $\D$.




For notational convenience, we set $X_{i',n} = {1 \over n} \sum_{i \leq i'} i(i-2) D_i'(n)$.
Clearly, $ Q_n' \geq X_{i',n}$.
On the other hand, (\ref{Assumption2}) implies that for every $\eps >0$ there exists $i_0 = i_0 (\eps)$ such that
whenever $i'>i_0$
for $n$ sufficiently large
\begin{equation} \label{tail}
{1\over n} \sum_{i > i'} i(i-2) D_i(n) < \eps.
\end{equation}
Since $D_i'(n) \leq D_i(n)$, for any $i'>i_0$ we have
$$ Q_n' \leq X_{i',n} + \eps.$$

We shall prove that for every such $i'$, $\mu - a.s.$
\begin{equation} \label{PartialLimit}
\lim_{n \rightarrow \infty} X_{i',n} = \sum_{i\leq i'} i(i-2)\la_i^{bond}=:Q_{i'}' .
\end{equation}
In turn, this will imply that for every $i'> i_0$
$$ Q_{i'}' \leq \liminf_{n \rightarrow \infty} Q_n' \leq \limsup_{n \rightarrow \infty} Q_n' \leq Q_{i'}' + \eps, \; \mu - a.s. $$
Now, letting $i' \rightarrow \infty$ yields
$$ Q' \leq \liminf_{n \rightarrow \infty} Q_n' \leq \limsup_{n \rightarrow \infty} Q_n' \leq Q' + \eps, \; \mu - a.s. $$
Since the choice of $\eps$ is arbitrary, we may eventually deduce that:
\begin{equation}\label{limit}
\lim_{n \rightarrow \infty} Q_n' = Q', \; \mu - a.s.
\end{equation}

So let us focus on proving (\ref{PartialLimit}). This will follow, if we show that
for every $\eps > 0$
\begin{equation} \label{ProbSum}
\sum_{n} \p \left[\left| X_{i',n} - Q_{i'}'\right| > \eps \right] <
\infty.
\end{equation}
(See for example Lemma 6.8 in~\cite{Petrov}.)
We will deduce the above inequality, proving that the summands are $o(1/n^3)$.

Thus, we continue with estimating
$$\p \left[\left| X_{i',n} - Q_{i'}'\right| > \eps \right], $$
for some fixed $\eps>0$.
Note that for $n$ sufficiently large
 $$ \left| \ex [X_{i',n}] - Q_{i'}' \right| \leq {\eps \over 2}.$$
Thus
$$\p \left[\left| X_{i',n} - Q_{i'}'\right| > \eps \right]
\leq \p \left[\left| X_{i',n} - \ex [X_{i',n}]\right| > {\eps \over 2} \right].$$
If the latter is realised, then there exists $i\leq i'$ for which
$$ {1\over n}\left|D_i'(n) - \ex [D_i'(n)] \right| > {\eps \over 2 \sum_{i\leq i'} i(i-2)}.$$
Therefore, setting $\eps' = {\eps \over 2 \sum_{i\leq i'} i(i-2)}$ we have
$$\p \left[\left| X_{i',n} - \ex [X_{i',n}]\right| > {\eps \over 2} \right]
\leq
\sum_{i\leq i'} \p \left[{1\over n} \left|D_i'(n) - \ex [D_i'(n)] \right| > \eps' \right].$$

We now show that each summand is $o(1/n^3)$.
To do so, we will condition on the size
of $C$. Recall that by Lemma \ref{uni-surv} conditional on $|C|=2k$,
the set $C$ is uniformly distributed among the $2k$-subsets of
$P(n)$. Moreover, if we replace one of the points in $C$ with
another one that does not belong to $C$, then $D_i'(n)$ can change by
at most 2. Therefore, applying (\ref{Conc}) we obtain
uniformly for any $k \in I$:
\begin{eqnarray*}
\p \left[ |D_i'(n) - \ex [D_i'(n)]|> \eps' n \ | \
|C|=2k \right] &\leq &
2 \exp \left(-\frac{ \eps'^2 n^2 }{4 k} \right)\\
&\leq & 2  \exp \left(-\frac{\eps'^2 n^2 }{4 (Mp+\ln n \sqrt{n})}
\right)
= o\left(\frac{1}{n^3}\right).
\end{eqnarray*}
Therefore, by (\ref{CConc})
\begin{eqnarray*}
\lefteqn{\p \left[|D_i'(n) - \ex [D_i'(n)]|> \eps' n \right]=} \\
& & \sum_{k \in I} \p \left[|D_i'(n) - \ex [D_i'(n)]|>\eps' n \ | \
|C|=2k\right] \ \p [|C|=2k] + o\left(\frac{1}{n^3}\right)
= o\left(\frac{1}{n^3}\right),
\end{eqnarray*}
for every $i$.

%


Now that we have proved (\ref{limit}), we are ready to conclude the proof
that $p_c^{bond}=\hat{p}_{bond}$.
Let $E\subseteq \D$ be the event over which $\lim_{n\rightarrow
\infty}Q_n'=Q'$; recall that $\mu (E)=1$.
Let $(\dd'(n))_{n \in \mathbb{Z}^+}\in E$.
If we condition on $\dd'(n)$ being the degree sequence on $G'(n)$,
then Lemma~\ref{uniform} implies that $G'(n)$ is the multigraph that
arises as the projection of a uniformly random perfect matching on
$P'(n)$ onto $V_n$.

If $p< \hat{p}_{bond}$ then $Q' < 0$.
For an arbitrary $\eps >0$, we
define $A(n)$ to be the set of multigraphs on $V_n$ whose largest
component has no more than $\varepsilon n$ vertices, for an arbitrary
$\varepsilon \in (0,1)$.
In this case, Theorem~\ref{key} implies that
$\lim_{n \rightarrow \infty}
\p [G'(n) \in A(n) \ | \ D(G'(n))=\dd' (n)]=1$,
for every $\eps >0$.

On the other hand, if $p > \hat{p}_{bond}$ then $Q' >0$.
Again by Theorem~\ref{key}, we deduce that there exists $\eps >0$ such
that for any $(\dd ' (n))_{n \in \mathbb{Z}^+} \in E$, if we
define $A(n)$ to be the set of multigraphs on $V_n$ whose largest
component has at least $\varepsilon n$ vertices, then $\lim_{n \rightarrow \infty}
\p [G'(n) \in A(n) \ | \ D(G'(n))=\dd' (n)]=1$.

However, in either case we want to know the limit of $\p [G' (n) \in A(n) ]$ as
$n \rightarrow \infty$ without conditioning on the degree sequence.
If $\omega \in \D$, then we let $\pi_n(\omega)$ denote the projection
of $\omega$ onto its $n$-th factor - recall that this is a degree
sequence on $V_n$.
Thus this probability can be expressed as follows:
\begin{eqnarray*}
 \p [ G'(n) \in A(n) ] &=& \sum_{ \dd'(n)  \in {\bf D}_n}
\p [ G'(n) \in A(n) \ | \ D( G'(n))=\dd' (n)]  \
\p [D(G'(n))=\dd ' (n) ]
 \nonumber \\
&=& \sum_{\dd '(n) \in {\bf D}_n}
\int_{\{ \omega \in \D \ : \
\pi_n (\omega )= \dd '(n) \}}
\p [G'(n) \in A(n) \ | \ D(G'(n))=\dd ' (n) ]
\ \mu( d\omega )
 \nonumber \\
&=& \int \p [G'(n) \in A(n) \ | \ D( G'(n))=
\pi_n(\omega) ] \
\mu( d \omega ).
\end{eqnarray*}
Since the integrand is bounded below by 0, applying Fatou's Lemma, we obtain:
\begin{eqnarray*}
\liminf_{n \rightarrow \infty}
\lefteqn{\p [ G'(n) \in A(n)  ] = \liminf_{n \rightarrow \infty}
\int \p [G'(n) \in A(n) \ | \ D( G'(n))=
\pi_n(\omega)] \ \mu( d \omega
 )} \nonumber \\
&\geq & \int\liminf_{n \rightarrow \infty}  \p [ G'(n) \in A(n)
\ | \ D( G'(n))= \pi_n(\omega) ]
\ \mu( d \omega )
\nonumber \\
&=& \int_{E} \liminf_{n \rightarrow \infty}
\p [ G'(n) \in A(n) \ | \ D( G'(n))=
\pi_n(\omega) ]
 \ \mu( d \omega)
= \int_{E}\mu( d \omega ) = 1.
\end{eqnarray*}
Now, applying the Reverse Fatou's Lemma (since the integrand is
bounded above by 1), we have
\begin{eqnarray*}
\lefteqn{\limsup_{n \rightarrow \infty}
\p [ G'(n) \in A(n) ] = \limsup_{n \rightarrow \infty}
\int \p [ G'(n) \in A(n) \ | \ D( G'(n))=
\pi_n(\omega )] \ \mu( d (\omega )) }
\nonumber \\
&\leq & \int\limsup_{n \rightarrow \infty}
\p [ G'(n) \in A(n) \ | \ D( G'(n))=
\pi_n(\omega) ] \
\mu( d \omega )
\nonumber \\
&=& \int_{E} \limsup_{n \rightarrow \infty}
\p [ G'(n) \in A(n) \ | \ D( G'(n))=
\pi_n(\omega)]
\ \mu( d \omega )
= \int_{E}\mu( d \omega ) = 1.
\end{eqnarray*}
The last two inequalities along with Lemma \ref{simple} complete the
proof of Theorem~\ref{main} for the bond percolation process.

\section{Site percolation}
In this section, we are dealing with site percolation, where for $p
\in (0,1)$ fixed, we make each vertex of $\tilde{G}(n)$ isolated
with probability $1-p$, independently of every other vertex,
deleting all of the edges that are attached to it. We will be
referring to this process as the {\em deletion} of the vertices.
This process applied to $\tilde{G}(n)$ induces a random degree
sequence on $V_n$, which, as in the previous section, we denote by
$\dd'(n)$. Now consider the effect of the deletion on the uniformly
random perfect matching on $P(n)$: if a vertex is deleted, then the
points of $P(n)$ that are the end-points of the edges attached to
that vertex are deleted (i.e. we remove them from $P(n)$).
Eventually, we are left with a set points of $P(n)$ that are the
end-points of the remaining edges, and we denote it by $C$.
Also, we let $P'(n)=P(\dd'(n))$.
As in the case of bond percolation, we establish a
bijection between the set of
perfect matchings between the points in $C$
and the set of perfect matchings on $P'(n)$. In turn, this gives rise to a bijection between
the set of perfect matchings on $C$ and the set of perfect matchings on $P'(n)$.

In the present setting, the set $\D_n$
consists of the degree sequences that are induced by the random
deletion of the vertices. The probability of a certain degree sequence in
$\D_n$ will therefore be the probability that this is the induced
degree sequence after the deletion.

We now argue that conditional on the choice of the points of $C$, each perfect
matching on $C$ has the same probability.
The perfect matching which is
realised after the deletion of the vertices is obtained in two
independent stages: firstly the uniform perfect matching on $P(n)$ is
realised and afterwards the random deletion of the vertices takes
place. It is the independence that allows us to consider these two
random experiments in reverse order. Thus, we choose first those
vertices that will be deleted and then we realise the perfect matching
on $P(n)$. Let $P_1(n)$ and $P_2(n)$ denote the sets of points
corresponding to the deleted vertices and the vertices that remain,
respectively.
Let $B$ be the subset of points in $P_2(n)$ that are matched with
points in $P_1(n)$.
Observe now that conditioning on the choice of $P_1(n)$ and $B$
is equivalent to conditioning on the choice of $C$, as the disjoint
union of $B$ and $C$ is $P_2(n)$.
Under this conditioning each perfect matching on $C$ has the same
probability, since the number of perfect matchings on the remaining points
is the same for every perfect matching on $C$.
Thus if $|C|=2k$, then each perfect matching on $C$ has probability $k!2^{k}/(2k)!$.

For a degree sequence $\dd'$, we let $S_{\dd'}$ be the set of subsets
of $P(n)$ which realise $\dd'$. Assume that the sum of the degrees in $\dd'$ is $2k$.
Note that if $\dd'(n) = \dd'$, then $P'(n)=P(\dd')$.
Hence if $m$ is a perfect matching on $P(\dd')$, then
\begin{eqnarray*}
\p [m \ | \ \dd'(n) =\dd'] &=& \sum_{C' \in S_{\dd'}} \p [m \ | \ C=C', \ \dd'(n) =\dd']~
\p [C=C' \ | \ \dd'(n)=\dd'] \\
&=& {k!2^k \over (2k)!} ~ \sum_{C' \in S_{\dd'}} \p [C=C' \ | \ \dd'(n)=\dd'] = {k!2^k \over (2k)!}.
\end{eqnarray*}

The parameter $p_c^{site}$ will be determined by $\dd' (n)$. If
$D_i'(n)$ denotes the number of vertices that have degree $i$ in $\dd '(n)$, then
letting
\begin{equation} \label{LimitEx}
\la_d^{site} = \lim_{n \rightarrow \infty}{1\over n} \ex [D_i'(n)]
\end{equation}
we shall prove that this limit exists for every $i\in \Nat$.
In fact we show that
\begin{equation} \label{Percs1}
\la_d^{site} = p \la_d^{bond}.
\end{equation}
This implies that
\begin{equation} \label{Percs}
\sum_i i(i-2)\la_i^{site} = p \sum_i i(i-2)\la_i^{bond}.
\end{equation}
Let $\hat{p}_{site}$ be the root of $\sum_i i(i-2)\la_i^{site}=0$. We will show that
$p_c^{site} = \hat{p}_{site}$. Then, by
(\ref{Percs}), we will deduce that $p_c^{site}=p_c^{bond}$.

We now prove the existence of the limit in (\ref{LimitEx}).
First of all we estimate the number of points in $P_2(n)$. Then we shall condition on 
a certain realisation of $P_2(n)$ and afterwards we
shall condition on the size of $B$ (i.e. on the size of $C$). From
this we will be able to estimate $\ex [D_i'(n)]$.
Let $D_d''=D_d''(n)$ be the number of vertices of degree $d$ that
survive the deletion - therefore $\ex [D_d'']=D_dp$.
The total degree in  $P_2(n)$ is
$M_2 = \sum_{d=1}^{\Delta}dD_d''$ and the linearity of expectation yields
 $\ex [M_2 ]=\sum_{d=1}^{\Delta}dD_d p=2Mp$.
As for every $n$ the maximum degree
in $\dd(n)$ is no more than $n^{1/9}$ a bounded differences inequality
(see for example Theorem 5.7 in~\cite{McD}) yields
\begin{equation} \label{M2Conc}
\p [|M_2 - \ex [M_2]| > n^{2/3}\ln n]\leq
\exp \left(- \Omega (\ln^2 n) \right).
\end{equation}

Now note that if we condition on $|B|=b$, then any $b$-subset of $P_2(n)$
is equilikely to occur as the set $B$ and it is the points of $B$ that
are deleted along with the points of $P_1(n)$.

Thus  if $P_{i-d}(i)$ denotes
the probability that after the random allocation of $B$ a certain
vertex in $P_2(n)$ of degree $i$ loses $i-d$ points, thus becoming a vertex of
degree $d$, the expected value of $D_d'$ is
\[ \ex [D_d'] = \sum_{i=d}^{\Delta} \ex[D_i''] P_{i-d}(i). \]
But for every $\eps >0$, there
exists $i_0$ such that $\sum_{i>i_0}D_i \leq \eps n$, for $n$
sufficiently large. Since $D_i''(n) \leq D_i(n)$ we obtain
\begin{equation}\label{PartialSums}
\sum_{i=d}^{i_0} \ex[D_i''] P_{i-d}(i) \leq
\ex [D_d'] \leq  \sum_{i=d}^{i_0} \ex[D_i''] P_{i-d}(i) + \eps n.
\end{equation}

We now calculate $P_{i-d}(i)$, for $d\leq i\leq i_0$.
Suppose that $M_2 =m_2$ and $|B|=b$.
Then if $P_{i-d}(i,b,m_2)$
is the conditional probability that after the random choice
of the $b$ points in $P_2(n)$, which has $m_2$ points,
a certain vertex of degree $i$ loses $i-d$ points we have
\begin{eqnarray*}
P_{i-d}(i,b,m_2) &=& {i \choose i-d} \frac{{ m_2 - i \choose b-i+d}}
{{m_2 \choose b}} =  {i \choose i-d} \ \frac{(m_2-i)!}{m_2!} \
\frac{b!}{(b-i+d)!} \ \frac{(m_2-b)!}{(m_2-b-d)!}.
\end{eqnarray*}
We shall assume that for any $n$ sufficiently large
$b \in [2Mp(1-p)-n^{2/3}\ln^2 n,2Mp(1-p)+ n^{2/3}\ln^2 n]$.
Indeed the following holds:
\begin{lemma} \label{edges-across}
Conditional on $M_2 \in I':=[2Mp- n^{2/3} \ln
n,2Mp+ n^{2/3}\ln n]$ we have
$b\in I:= [2Mp(1-p)- n^{2/3} \ln^2 n ,2Mp(1-p)+ n^{2/3} \ln^2 n ]$ with
probability $1-\exp \left(-\Omega(\ln^2 n) \right)$.
\end{lemma}
\begin{proof}
Assume that $M_2 = m_2$ for some $m_2\in I'$. Therefore, $M_1 = 2M -
m_2 \in [2M(1-p)- n^{2/3} \ln n,2M(1-p)+ n^{2/3}\ln n]$. We shall also
condition on a particular realisation of the sets $P_1(n)$ and $P_2(n)$.

The probability that a certain point in $P_2(n)$ is adjacent to a
point in $P_1(n)$ is ${M_1 \over 2M-1 }= (1-p)(1+O(n^{-1/3} \ln n))$.
Then $\ex [b]=2Mp(1-p)\pm O(n^{2/3}\ln n)$.

We now show that $b$ is concentrated around its expected value, using
Theorem 7.1 from~\cite{McD}. We first describe here the more general
setting on which this theorem applies and afterwards we will consider
$b$.

Let $W$ be a finite probability space that is also a metric space with
its metric denoted by $d$. Assume that $P_0,\ldots ,P_s$ is
 a sequence of partitions on $W$, such
that $P_{i+1}$ is  a refinement of $P_i$, $P_0$ is the partition
consisting of only one part, that is $W$, and $P_s$ is the partition
where each part is an element of $W$. Assume that whenever
$A, B \in P_{i+1}$ and $C\in P_i$
are such that $A,B \subseteq C$, then there is a bijection $\phi: A
\rightarrow B$ such that $d(x,\phi (x))\leq c_i$.

Now, let $V$ be a uniformly random element of $W$ and let
$f:W\rightarrow \mathbb{R}$ be a function on $W$ satisfying
$|f(x)-f(y)|\leq d(x,y)$. Then
\begin{equation}\label{ConcM}
\p [|f(V)- \ex[f(V)]|> t] \leq 2 \exp\left(-2 {t^2\over \sum_{i=1}^s c_i^2
}\right).
\end{equation}

In our context the uniform space of all perfect matchings on $P(n)$
will play the role of $W$. Let $\textbf{M}$ denote it. Its metric
will be the symmetric difference of any two perfect matchings,
regarded as sets. It is easy to see that this satisfies the
properties a metric has by its definition. We shall consider a
series of partitions on $\textbf{M}$ denoted by $P_0,\ldots,
P_{M-1}$, where $P_0$ is $\textbf{M}$ itself and each part of
$P_{M-1}$ will be a perfect matching in $\textbf{M}$. To define the
$i$-th partition, we define an ordering on the edges of each perfect
matching. Consider first a linear ordering of all the points in
$P(n)$. This induces a linear ordering on the edges of a perfect
matching: if $e_1$ and $e_2$ are two edges, then $e_1 <e_2$ if the
smallest point in $e_1$ is smaller than the smallest point in $e_2$.
Now a part of $P_i$ consists of those perfect matchings whose $i$
smallest edges are a particular set of $i$ edges, provided that such
a set of perfect matchings is non-empty. We call such an $i$-set of
edges a \emph{prefix}. Moreover, given a perfect matching, we call
its $i$ smallest edges its $i$-\emph{prefix}.

Given such an $i$-set of edges, let $C$ be the set of perfect
matchings that have these $i$ edges as their $i$-prefix.
Now consider two $i+1$-subsets that contain this
$i$-set and are both prefixes. Suppose that $e_A$ and $e_B$ are the last
edges on which they differ.
Let $A$ and $B$ respectively denote the sets
of perfect matchings that have these two $i+1$-sets as their
$i+1$-prefixes.

There is a natural bijection $\phi:A \rightarrow B$ between them.
Observe first that the
smallest vertex in $e_A$ and $e_B$ is the same. In particular, let us assume that
$e_A=(x,y_A)$ and $e_B=(x,y_B)$. If $m$ is a matching in $A$, then
$\phi (m)$ is the matching in $B$, where $y_A$ is adjacent to the
vertex that $y_B$ was adjacent to in $m$; every other edge remains unchanged.
Note that the symmetric difference of $m$ and $\phi(m)$ is 4. In other
words, $c_i=4$.

Now, we are ready to apply the concentration bound (\ref{ConcM}) to
$b$. For any $m \in \textbf{M}$, we let
$b(m)$ be the number of edges between $P_1(n)$ and $P_2(n)$. Observe
that for any two perfect matchings $m,m' \in \textbf{M}$, always
$|b(m)-b(m')|$ is no more than the size of the symmetric difference of
$m$ and $m'$. Thus applying (\ref{ConcM}) with $t= n^{2/3}\ln^2 n /2$, the
lemma follows, for $n$ large enough.
\end{proof}
Thus, uniformly for any $b \in I$ and any $m_2 \in I'$
we have:
\begin{eqnarray*}
P_{i-d}(i,b, m_2) &=& {i \choose i-d} \frac{b^{i-d}(m_2-b)^d}{m_2^i}
\left(1+O\left(\frac{1}{n}\right) \right) \\
&=& {i \choose d} (1-p)^{i-d} p^d
\left(1+O\left(\frac{\ln^2 n}{n^{1/3}}
\right) \right).\end{eqnarray*}
A standard concentration argument shows that uniformly for any $i\leq i_0$
we have 
\begin{equation}\label{DiConc} \p [|D_i''(n) - \ex[D_i''(n)]|\geq \sqrt{n}\ln n] \leq
\exp(-\Omega (\ln^2 n)). \end{equation} 
Thus if we also  set
$I''(i)=[\max\{pD_i-\ln n\sqrt{n}, 0\},pD_i+\ln n\sqrt{n}]$, we have
\begin{equation}\label{Intervals}
\p [b \not \in I \ \mbox{or} \ M_2 \not \in I' \mbox{or}\ D_i'' \not
\in I''(i),\ \mbox{for some $i\leq i_0$}]=o(n^{-3}).\end{equation}
Therefore, the right-hand side of (\ref{PartialSums}) becomes
\begin{eqnarray*}
\ex [D_d'] &\leq& \sum_{k \in I} \sum_{k' \in I'}\sum_{i=d}^{i_0}
\sum_{k_i'' \in I''(i)}
 k_i'' P_{i-d}(i,k,k') \p [b=k, \ M_2=k', \ D_i''=k_i''] + \eps n
+o(n^{-2})  \\
&=& \sum_{i=d}^{i_0}\sum_{k_i'' \in I''(i)} k_i''
{i \choose d} (1-p)^{i-d} p^d \p [D_i''=k_i'']
\left(1+O\left(\frac{\ln^2 n}{n^{1/3}} \right) \right)
+\eps n + o(n^{-2}).
\end{eqnarray*}
But by (\ref{DiConc}), we have $\sum_{k_i'' \in I''(i)} k_i'' \p [D_i''=k_i''] =\ex
[D_i''(n)]-o(n^{-2})=D_i(n)p-o(n^{-2})$.
Substituting this into the above expression, (\ref{PartialSums}) now
yields:
$$
\ex [D_d'] \leq  p \sum_{i=d}^{i_0}D_i(n)
{i \choose d} (1-p)^{i-d} p^d
\left(1+O\left(\frac{\ln^2 n}{n^{1/3}} \right) \right)
+\eps n + o(n^{-2})
$$
and also repeating the above estimations,
$$
\ex [D_d'] \geq  p \sum_{i=d}^{i_0}D_i(n)
{i \choose d} (1-p)^{i-d} p^d
\left(1+O\left(\frac{\ln^2 n}{n^{1/3}} \right) \right)
+ o(n^{-2}).
$$
Therefore,
$$p \sum_{i=d}^{i_0}\la_i
{i \choose d} (1-p)^{i-d} p^d
\leq \liminf_{n \rightarrow \infty}{1\over n} \ex [D_d']$$
and
$$ \limsup_{n \rightarrow \infty}{1\over n} \ex [D_d'] \leq
p \sum_{i=d}^{i_0}\la_i
{i \choose d} (1-p)^{i-d} p^d +\eps. $$
Letting $i_0\rightarrow \infty$ and $\eps \rightarrow 0$ we obtain
\begin{equation} \label{LimitFinal}
\la_{d}^{site}\equiv \lim_{n \rightarrow \infty}{1\over n} \ex [D_d']=p
\sum_{i=d}^{\infty}\la_i {i \choose d} (1-p)^{i-d} p^d,
\end{equation}
which yields (\ref{Percs1}) through (\ref{labond}).

Now we let $Q_n' = {1\over n}\sum_{i\geq 2} i(i-2)D_i'(n)$. We will
show that $\mu - a.s.$
\begin{equation}\label{LimitSite}
\lim_{n \rightarrow \infty} Q_n' = \sum_{i\geq
1}i(i-2)\la_i^{site}\equiv Q'.
\end{equation}

To prove this we argue as in the case of bond percolation: setting
$X_{i',n}={1\over n}\sum_{i\leq i'}i(i-2)D_i'(n)$, for every $\eps >0$
and any $i'$ large enough we have
$$X_{i',n}\leq Q_n' \leq X_{i',n}+\eps, $$
if $n$ is also large enough. (Obviously the first inequality holds for
every $i'$ and $n$.)
Thus the existence of the $\mu -a.s.$ limit of $Q_n'$ will be
established once we show that for any $i'$ $\mu - a.s.$
$\lim_{n\rightarrow \infty}X_{i',n}=\sum_{i\leq i'}
i(i-2)\la_i^{site}\equiv Q_{i'}'$. We then let $i'\rightarrow
\infty$ and $\eps \rightarrow 0$ to deduce that $\mu - a.s.$ $\lim_{n
\rightarrow \infty} Q_n' = Q'$.

The almost sure convergence of $X_{i',n}$ to $Q_{i'}'$ can be shown as
in the case of bond percolation. In other words, we need to prove that
the condition in (\ref{ProbSum}) is satisfied in the present context.
As before, we will show that for every
$i\leq i'$ the random variable $D_i'(n)$ is sharply concentrated around its expected
value: that is its tails converge to 0 exponentially fast.
Recall that the total degree in $P_2(n)$ is denoted by $M_2$.

Conditional on a certain realisation of $P_2(n)$,
with $M_2 = m_2$ for some $m_2 \in I'$, and $|B|=b$ for some $b \in I$,
the value of $D_i'(n)$ is determined by the random choice of the set $B$
in $P_2(n)$.
Note that $D_i'$ can change by at most 2, if we replace one
element of $B$ by another one. Therefore we may apply (\ref{Conc}):
\begin{eqnarray*}
\lefteqn{\p \left[|D_i'(n) - \ex [D_i'(n)]|> \ln n \sqrt{n} \ | \
|B|=b, \ P_2(n),\ |P_2(n)|=m_2\right]\leq}\\
& &4\exp \left( -\frac{n \ln^2 n}{2 (2Mp(1-p)+n^{2/3}\ln^2 n)}\right) =
\exp \left( - \Omega ( \ln^2 n) \right),
\end{eqnarray*}
uniformly for any $b \in I$ and $m_2 \in I'$.
Hence, the above inequality along with (\ref{M2Conc}) and Lemma (\ref{edges-across}) imply that
\begin{eqnarray*}
\p \left[|D_i'(n) - \ex [D_i'(n)]|> \ln n \sqrt{n} \right] &=& o(n^{-3}).
\end{eqnarray*}
Since $i\leq i'$ and $i'$ is bounded, condition (\ref{ProbSum}) is satisfied, and
therefore,  $\mu - a.s.$
$\lim_{n\rightarrow \infty}X_{i',n}=\sum_{i\leq i'}
i(i-2)\la_i^{site}\equiv Q_{i'}'$. Now this concludes the proof of
(\ref{LimitSite}).

The proof that $p_c^{site}=\hat{p}_{site}$ 
is identical to that for $\hat{p}_{bond}$, and it is omitted.

\end{document}